%


\input mssymb


%

%

\newcount\skewfactor
\def\mathunderaccent#1#2 {\let\theaccent#1\skewfactor#2
\mathpalette\putaccentunder}
\def\putaccentunder#1#2{\oalign{$#1#2$\crcr\hidewidth
\vbox to.2ex{\hbox{$#1\skew\skewfactor\theaccent{}$}\vss}\hidewidth}}



\def\rest{\mathord{\restriction}}

\def\phi{\varphi}

\def\su{\subseteq}
\def\a{\alpha}

\def\k{\kappa}

\def\om{\omega}

\def\lng{\langle}
\def\rng{\rangle}

\def\sm{\setminus}
\def\cont{{2^{\aleph_0}}}


\def\deg{{\rm deg}}


\def\endproof#1{\hfill  
{\parfillskip0pt$\smiley_{\hbox{{#1}}}$\par\medbreak}}

\def\proof{\smallbreak\noindent{\sl Proof}: }


\def\Cal#1{{\cal #1}}
\def\isom{\cong}


\newbox\noforkbox \newdimen\forklinewidth
\forklinewidth=0.3pt   

\setbox0\hbox{$\textstyle\bigcup$}
\setbox1\hbox to \wd0{\hfil\vrule width \forklinewidth depth \dp0
			height \ht0 \hfil}
\wd1=0 cm
\setbox\noforkbox\hbox{\box1\box0\relax}
\def\unionstick{\mathop{\copy\noforkbox}\limits}
\def\nonfork#1#2_#3{#1\unionstick_{\textstyle #3}#2}
\def\nonforkin#1#2_#3^#4{#1\unionstick_{\textstyle #3}^{\textstyle  
#4}#2}

\setbox0\hbox{$\textstyle\bigcup$}
\setbox1\hbox to \wd0{\hfil$\nmid$\hfil}
\setbox2\hbox to \wd0{\hfil\vrule height \ht0 depth \dp0 width
				\forklinewidth\hfil}
\wd1=0cm
\wd2=0cm
\newbox\doesforkbox
\setbox\doesforkbox\hbox{\box1\box0\relax}
\def\nunionstick{\mathop{\copy\doesforkbox}\limits}

\def\fork#1#2_#3{#1\nunionstick_{\textstyle #3}#2}
\def\forkin#1#2_#3^#4{#1\nunionstick_{\textstyle #3}^{\textstyle  
#4}#2}

\font\circle=lcircle10

\setbox0=\hbox{~~~~~}
\setbox1=\hbox to \wd0{\hfill$\scriptstyle\smile$\hfill} 
\setbox2=\hbox to \wd0{\hfill$\cdot\,\cdot$\hfill} 

\setbox3=\hbox to \wd0{\hfill\hskip4.8pt\circle i\hskip-4.8pt\hfill}  


\wd1=0cm
\wd2=0cm
\wd3=0cm
\wd4=0cm

\newbox\smilebox
\setbox\smilebox \hbox {\lower 0.4ex\box1
		 \raise 0.3ex\box2
		 \raise 0.5ex\box3
		\box4
		\box0{}}
\def\smiley{\leavevmode\copy\smilebox}

\headline={\tenrm  
\number\folio\hfill\jobname\hfill\number\day.\number\month.\number\year}


\outer\long\def\ignore#1\endignore{}

\newcount\itemno
\def\itm{\advance\itemno1 \item{(\number\itemno)}}
\def\ritm{\advance\itemno1 \item{)\number\itemno(}}

\def\aitm{\advance\itemno1 

\item{(\letter\itemno)}}
^^L

\def\letter#1{\ifcase#1 \or a\or b\or c\or d\or e\or f\or g\or h\or
i\or j\or k\or l\or m\or n\or o\or p\or q\or r\or s\or t\or u\or v\or
w\or x\or y\or z\else\toomanyconditions\fi}
\def\raitm{\advance\itemno1 \item{)\rletter\itemno(}}
\def\rletter#1{\ifcase#1\or `\or a\or b\or c\or d\or e\or f\or g\or
h\or i\or k\or k\or l\or n\or n\or q\or r\or t\or v
\else\toomanyconditions\fi}


\newcount\secno
\newcount\theono

\catcode`@=11
\newwrite\mgfile

\openin\mgfile \jobname.mg
\ifeof\mgfile \message{No file \jobname.mg}
	\else\closein\mgfile\relax\input \jobname.mg\fi
\relax
\openout\mgfile=\jobname.mg

\newif\ifproofmode
\proofmodetrue            

\def\@nofirst#1{}

\def\neusection{\advance\secno by 1\relax \theono=0\relax}
\def\neuchap{\secno=0\relax\theono=0\relax}

\neuchap

\def\labelit#1{\global\advance\theono by 1%
             \global\edef#1{%
             \number\secno.\number\theono}%
             \write\mgfile{\@definition{#1}}%
}

^^L




\def\ppro#1#2:{%
\labelit{#1}%
\smallbreak\noindent%
\@markit{#1}%
{\bf\ignorespaces {#2}:}}





\def\@definition#1{\string\def\string#1{#1}
\expandafter\@nofirst\string\%
(\the\pageno)}

\def\@markit#1{
\ifproofmode\llap{{ \expandafter\@nofirst\string#1\ }}\fi%
{\bf #1\ }
}

\def\h@markit#1{
\ifproofmode\edef\nxt{\string#1\ }%
{\tenrm\beginL\nxt\endL}
\fi%
{{\bf\beginL #1\endL}}
}
 ^^L
\def\labelcomment#1{\write\mgfile{\expandafter
		\@nofirst\string\%---#1}}

\catcode`@=12

\def\refrence#1#2:{\write\mgfile{\def\noexpand#1{#2}}%
\areference{#2}}
\def\areference#1{\medskip\item{[#1]} \ignorespaces}

\newcount\referencescount
\def\numrefrence#1{\advance\referencescount1
\edef#1{\number\referencescount}%
\write\mgfile{\def\noexpand#1{#1}}%
\areference{#1}}

\def\numericalreferences{\let\refrence\numrefrence}

\newcount\scratchregister
\def\simplepro{\scratchregister\theono\advance\scratchregister1 

\edef\scratchmacro{\number\secno.\number\scratchregister}%
\expandafter\ppro\csname\scratchmacro\endcsname}


\font\teneuf=eufm10
\font\seveneuf=eufm7
\font\fiveeuf=eufm5
\newfam\euffam
\textfont\euffam=\teneuf
\scriptfont\euffam=\seveneuf
\scriptscriptfont\euffam=\fiveeuf


\def\edge{\vrule width0.7cm height3pt depth-2.6pt}



\def\complete#1#2#3{\vtop{\offinterlineskip
\ialign{
&\hfil$##$\hfil\cr
\,#1\,&\omit\hfil\edge\hfil&\,#2\,\cr
\noalign{\vskip2pt}
 {\Big|} & #3  & {\Big|}\cr
\omit\hfill\hrulefill&\omit\hrulefill&\omit\hrulefill\hfill\cr
&  \qquad \cr}}}

\def\highway#1#2#3{\hbox{$\,#1\,$\edge$\,#2\,$\edge$\,\cdots\,$\edge$ 
\,#3\,$}}

%
\def\bridge#1#2#3#4#5#6#7{%
\vcenter{\ialign{&\hfil$##$\cr
#1\cr
&\setminus\cr
\cr
&/\cr
#2\cr
}}\vcenter{\highway{#3}{#4}{#5}}
\vcenter{\ialign{&\hfil$##$\cr
&#6\cr
/\cr
\cr
\setminus\cr
&#7\cr
}}}



\def\dead#1#2#3#4#5#6{\highway{#1}{#2}{#3}\edge\complete{#4}{#5}{#6}}










\magnification=1200
\baselineskip18pt
\proofmodefalse

\font\bigfont cmbx10 scaled \magstep2
\font\namefont cmbx10 scaled \magstep2
{{
\obeylines

\everypar={\hskip0cm plus 1 fil}
\bf
\parskip=0.4cm

{\bigfont UNIVERSAL BRIDGE FREE GRAPHS}

\bigskip
{\rm July 1994}
\bigskip

\parskip0.1cm

{\namefont Martin Goldstern}
Institute of Algebra
Technical University
Vienna, Austria
{\tt goldstrn@email.tuwien.ac.at}

\bigskip

{\namefont Menachem Kojman}
Department of Mathematics 

Carnegie Mellon University
Pittsburgh, PA 15213, USA
{\tt kojman@andrew.cmu.edu}

}
\vfill
\everypar{}
\rm
\leftskip2cm\rightskip2cm

ABSTRACT. We prove that there is no countable universal $B_n$-free
graph for all $n$ and that there is no countable universal graph in
the class of graphs omitting all cycles of length at most $2k$ for
$k\ge 2$.

\footline{\hfill}
\global\pageno0

\eject
}

\noindent
{\bf \S0\ Introduction}

Several papers have addressed the problem of existence of a universal
element among all countable graphs omitting given finite subgraphs
(see [KP] and the comprehensive bibliography there, and the most
recent [CK] and [KP1]).

Given a graph $F$ we say that a graph $G$ is $F$-free if $F$ is not
isomorphic to a subgraph of $G$. A countable $F$-free graph $G^*$ is
universal (strongly universal) in the class of all countable $F$-free
graphs if every countable $F$-free graph is isomorphic to a subgraph
(an induced subgraph) of $G^*$.

In [CK], Cherlin and Komjath raise the problem of determining for  
which
finite trees $T$  there exists a universal countable $T$-free
graphs. In this paper we describe an infinite set of finite trees
$B_n$, which we call bridges, and show that for no $n$ is there a
universal countable  $B_n$-free graph.

In [CK] it is proved that for all $n\ge 4$ there is no universal
countable $C_n$-free graph (where $C_n$ is a cycle of length $n$). In
[KMP] it is proved, on the other hand, that a strongly universal
countable graph exists among all countable graphs that omit all odd
cycles of length at most $2k+1$. What if we intersect some of those
classes, say look at all graphs omitting $C_3,C_4,C_5,C_6$? We show
here, using an idea of S.~Mozes, that when all cycles of
length at most $2k$ are to be excluded (for $k\ge 2$), then there is
no universal countable graph.

\ppro \problem Problem: Is there a countable universal graph in the
class of graphs omitting all cycles of length at most $2k+1$ (for
$k\ge 2$)? More generally, for what sets $F\su N$ does the
class of graphs omitting $\{C_n:n\in F\}$ have a countable universal
element?

Following [KP] we make the following definition:

\ppro \complexity Definition: Let $\Cal G$ be a class of graphs. 

\item{(i)} The
complexity ${{\rm cp}}({\Cal G})$ of $\Cal G$ is the minimal
cardinality $\k$ of a set $I$ of graphs in $\Cal G$ with the property
that every member in $\Cal G$ is embedded as an induced subgraph into
at least one of the members of $I$.

\item{(ii)} The weak complexity ${{\rm wcp}}(\Cal G)$ is defined by
omitting the word ``induced'' from the definition of ${\rm cp}(\Cal
G)$.

\medbreak
\noindent
{\bf Notation}: We denote the vertex degree of $v$ in a graph $G$ by
$\deg_G(v)$. The {\it length} of a path is the number of edges in the
path.  For all $m$ let $K_m$ denote a complete graph with $m$  
vertices
$k_1,\ldots,k_m$, let $P_m$ be a simple path of length $m$ with
vertices $p_0,\ldots ,p_m$ and edges $(p_i,p_{i+1})$ for $i\le m$ and
let $C_m$ denote a cycle of length $m$. Let us call a simple path
$h_0,h_1\ldots,h_k$ in a graph $G$ a {\it highway} iff  
$\deg_G(v_i)=2$
for all $0<i<k$.

\ppro \advice Advice: Drive carefully.

\neusection 

\bigbreak
\noindent
{\bf \S\number\secno\ Bridge-free graphs}

\ppro \bridgeDef Definition: A finite graph with $n+5$ vertices
is called an $n$-bridge iff it is isomorphic to $B_n=\lng V,E\rng$
where $V=\{a,b,c,x_1,x_2,\ldots,x_n,d,e\}$ and $E=\{(a,c),(b,c),
(x_n,d),x_n,e)\} \cup \{(x_i,x_{i+1}) : 1 \le i<n\}$

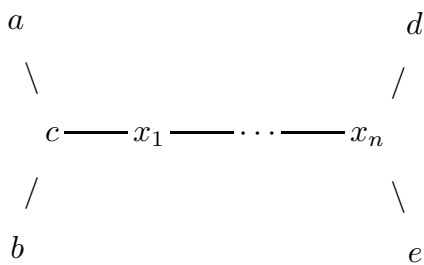
$$     \bridge abc{x_1}{x_n}de   $$

\medbreak

\ppro \deadEnd Definition: Let us call a graph $D$ a {\it dead end}  
if it is
isomorphic to $K_{n+3}$ to which a simple path $P_{n+1}$ is freely
adjoined by identifying $k_{n+3}$ with $p_0$.

$$ \dead{p_{n+1}}{ }{p_1}{k_{n+3}}{k_{1}}{K_{n+3}} $$

\medbreak

\ppro \rigidOne Claim: A dead end is $B_n$-free and if a dead end $D$
with vertices $k_1,\ldots ,k_{n+3}=p_0,x_1,\ldots p_{n+1}$ is a
subgraph of a $B_n$-free graph $G$ then $\deg_D(v)=\deg_G(v)$ for all
vertices $v\in D$ except maybe $v=p_{n+1}$.

\proof:  Suppose first that for some $i\le n$ there is an edge
$(p_i,y)$ in $G$ which is not an edge of $D$. If $y\notin D$, by
labeling $y$ as $a$, labeling $p_i$ as $c$ and $p_{i+1}$ as $b$ it is  
possible
to label vertices of $D$ as $x_i$ ($i\le n$) and as $d,e$ to produce  
a
copy of $B_n$.

Suppose, then, that $(p_i,p_j)$ is an edge not among the edges of  
$D$.
Without loss of generality, $j\ge i+2$. Now label $p_j$ as $a$, label
$P_i$ as $c$ and $p_{i+1}$ as $b$. Again, a copy of $B_n$ is easily
found.

The remaining possibility for an edge of the form $(p_i,y)$ is that
$y=k_j$ for some $j\le n+3$.  Here we distinguish two
subcases. First, $j=1$ (and, of course, $i>1$). Labeling $p_1$ as  
$a$,
$k_1$ as $c$ and $p_1$ as $b$, the remaining $n+2$ vertices of
$K_{n+3}$ complete the three labeled ones to make a copy of $B_n$.

Second, $j\not=1$. In this case if $i>1$ label $k_j$
as $c$, label $k_1$ as $x_n$ and label $p_1$ as $e$. The remaining
vertices of $K_{n+3}$ serve as $x_i$ for $1\le i < n$ and as $b$. If,
however, $i=1$, label $p_1$ as $c$, label $p_2$ as $a$ and label  
$k_j$
as $b$. Again, a copy of $B_n$ is found.

We show next that the $\deg_D(k_i)$ is preserved.  Suppose that
$(k_i,y)$ is an edge in $G$ which is not an edge in $D$.  We already
proved that $y\not=p_j$ for all $j\le n$. Therefore, either $y\notin  
D$
or $y=p_{n+1}$. If $i=1$ label $p_1$ as $a$, label $y$ as $b$ and
$k_1$ as $c$; otherwise label $k_i$ as $c$, $y$ as $a$, $p_1$ as $e$
and $k_1$ as $x_n$. In both cases a copy of $B_n$ results
\endproof\rigidOne

\ppro \driveThru Definition: Let us call a graph $T$ a
{\it drive through} if it is isomorphic to the graph obtained as
follows: Let $k_1,\ldots,k_{n+2}$ be the vertices of a copy of
$K_{n+2}$. For $1<i<n+2$ adjoin freely to $k_i$ a copy of a dead end
by identifying $p_{n+1}$ in that copy with $k_i$. To $k_1$ connect a
vertex $l$ by an edge and to $k_{n+2}$ connect a vertex $r$ by an
edge. Call $l$ the {\it left exit} of $T$ and call $r$ the {\it right
exit} of $T$.

PICTURE

\ppro \rigidTwo Claim: A  drive through $T$ is $B_n$-free and if $T$
is a subgraph of a $B_n$-free graph $G$ then $\deg_T(v)=\deg_G(v)$  
for
all vertices $v\in T$ except $l$ and $r$.

\proof: Suppose to the contrary that $B_n$ is a subgraph of a drive
through $T$. As $\deg_{B_n}(c)=\deg_{B_n}(x_n)=3$, both $a$ and $x_n$
are either in a copy of $K_{n+3}$ or in the copy of $K_{n+2}$. Both
cannot be in the same copy of $K_{n+3}$ because the minimum of
distances of $x_i$ to $c$ and $x_n$ is smaller than $n+1$, and all  
the
points satisfying this would be  in the same dead end as $x_n$ and  
$c$,
contrary to claim \rigidOne. Similarly, $c$ and $x_n$ are not both in
the copy of $K_{n+2}$.

Also, $c$ and $x_n$ cannot be in different copies of $K_{n+3}$, or in  
a
copy of $K_{n+3}$ and in the copy of $K_{n+2}$ because the distance
between $x_n$ and $c$  would be greater than $n$. We conclude that  
$T$ is
$B_n$-free.

Suppose that $T$ is a subgraph of a $B_n$ free graph $G$. By claim
\rigidOne\ we know that $\deg_T(v)=\deg_G(v)$ for all vertices $v$ in
the dead ends except those which are also in the copy of $K_{n+2}$.
Suppose that for some vertex $k_i$ in the copy of $K_{n+2}$ there is
an edge $(k_i,y)$ in $G$ which is not an edge of $T$. Label $y$ as
$a$. If $i=1$ or $i=n+2$ label $l$ or $r$ respectively as $b$. Label
$k_i$ as $c$. Label $n$ of the remaining $k_i$ as $x_1,\ldots,x_n$.
Label the last remaining $k_i$ as $d$. If this $i$ is $1$ or $n+1$
label $l$ or $r$ respectively as $e$.  Otherwise label as $e$ the
vertex $p_n$ in the dead end adjoined to $k_i$. This yields a copy of
$B_n$.\endproof\rigidTwo

 For every $\epsilon\in {}^{\om}2$ we construct a 

connected $B_n$-free graph $G_\epsilon$ as
follows.

Let $T^\epsilon(m)$ for $m\in N$ and $\epsilon\in{}^\om2$ be disjoint
copies of a drive through. Let $l^\epsilon(m)$ and $r^\epsilon(m)$ be
the left and right exists of $T^\epsilon(m)$. Let $D^\epsilon$ be a
copy of $D$ with vertices
$k_1^\epsilon,\ldots,k_{n+3}^\epsilon=p_0^\epsilon,\ldots
p_{n+1}^\epsilon$. Let $H^\epsilon(m)$ be a simple path of length
$n-1+\epsilon(m)$ with vertices
$h^\epsilon_0(m),\ldots,h^\epsilon_{n-1+\epsilon(m)}(m)$.

Adjoin $D^\epsilon$ to the $l^\epsilon(0)$ by setting
$p^\epsilon_{n+1}=l^\epsilon(0)$.  Connect $r^\epsilon(m)$ to
$l^\epsilon(M+1)$ by $H^\epsilon(m)$ by setting
$r^\epsilon(m)=h^\epsilon(0)$ and
$l^\epsilon(+1)=h^\epsilon_{n-1+\epsilon(m)}(m)$.  (If $n=1$ then  
when
$\epsilon(m)=0$ we identify $r^\epsilon(m)$ with $l^\epsilon(m+1)$.)

Let $G_\epsilon=D^\epsilon\cup \bigcup T^\epsilon(m)\cup \bigcup
H^\epsilon(m) $.

Let us observe that all highways in $G_\epsilon$ are either of length
$n+1$ or of length $n+2$. All highways that have an end of degree  
$n+3$
are of length $n+1$ except a unique highway --- the one containing
$l^\epsilon(0)$ --- which is of length $n+2$. Let us denote this
highway by $H(\epsilon)$.

\ppro \GeeEpsilon Claim: The graph $G_\epsilon$ is $B_n$-free and if
$G_\epsilon$ is a subgraph of a $B_n$-free graph $G$ then the vertex
degree of every vertex $v\in G_\epsilon$ in $G_\epsilon$ equals the
degree of $v$ in $G$.

\proof: A similar argument to that in \rigidTwo\ shows that  
$G_\epsilon$
is $B_n$-free. Suppose now that $G_\epsilon\su G$ and that $G$ is
$B_n$-free. By \rigidTwo\ we already know for
$\deg_{G_\epsilon}(v)=\deg_G(v)$ for each vertex $v\in T^\epsilon(m)$
except $l^\epsilon(m),r^\epsilon(m)$. If, however,
$\deg_G(v)>\deg_{G_\epsilon}(v)$ when $v$ is on one of the highways  
of
$G_{\epsilon}$, there must be some $y\in G\sm G_\epsilon$ such that
$(v,y)$ is an edge of $G$ and a copy of $B_n$ is easily
produced.\endproof\GeeEpsilon

\ppro \isom Corollary: For every $\epsilon\in {}^\om2$  and every  
connected
$B_n$-free graph $G$, if $G_\epsilon\su G$ then $G_\epsilon=G$.

\proof: Suppose that $y\in G\sm G_\epsilon$. By connectedness of $G$
we may assume that $y$ is connected by an edge to a vertex of
$G_\epsilon$. This contradicts \GeeEpsilon

\ppro \NonIsom Claim: If $\epsilon\not=\nu$ are two members of  
${}^\om
2$ then $G_\epsilon$ and $G_\nu$ are not isomorphic.

\proof: Suppose that $f:G_\epsilon\to G_\nu$ is an isomorphism. We
show that $\epsilon=\nu$. Clearly,  $f$ maps every highway in
$G_\epsilon$ onto some highway in $G_\nu$.

The highway $H(\epsilon)$ has to be mapped by $f$ onto $H(\nu)$, both
being the unique highways in their respective graphs of length $n+2$
with an end of degree $n+3$. As $l^\epsilon(0)$ is connected by an  
edge to
the end of $H(\epsilon)$ that has degree $n+2$, we conclude that
$f(l^\epsilon(0))=l^\nu(0)$. We argue by induction on $m$ that
$H^\epsilon(m)$ is mapped by $f$ onto the $H^\nu(m)$ and that
$f(l^\epsilon(m+1))=l^\nu(m+1)$.

 If $m=0$, we already showed that $f(l^\epsilon(0))=l^\nu(0)$.
Therefore $f(r^\epsilon(0))\not=l^\nu(0)$. Also, $f(r^\epsilon(0)$
cannot lie on any of the highways in $T^\nu(0)$ which are part of a
dead end, because both ends of $H^\epsilon(0)$ have degree $n+2$.
Therefore necessarily $f(r^\epsilon(0))=r^\nu(0)$ and consequently
$H^\epsilon(0)$ is mapped by $f$ onto the $H^\nu(0)$, with
$f(l^\epsilon(1))=l^\nu(1)$.

Similarly, if $f$ maps $H^\epsilon(m)$ onto $h^\nu(m)$ with
$f(l^\epsilon(m+1))=l^\nu(m+1)$, it follows that $f$ maps
$H^\epsilon(m+1)$ onto $h^\nu(m+1)$ with  
$f(l^\epsilon(m+2))=l^\nu(m+2)$.

As for all $m$ we have established that $n-1+\epsilon(m)=n-1+\nu(m)$,
we have shown that $\epsilon=\nu$.\endproof\NonIsom

\ppro \main Theorem: There is no universal $B_n$-free graph. In
fact, the weak complexity of the class of countable $B_n$-free graphs
equals $\cont$.

\proof
Suppose that $\{G_\a:\a\in I\}$ is a collection of less than $\cont$
many countable $B_n$ free graphs. By splitting each
graph to its connected components we assume that each $G_\a$ is  
connected.
Suppose that for every $\epsilon\in {}^\om2$ the graph $G_\epsilon$
constructed above is isomorphic to a subgraph of $G_\a$ for some
$\a\in I$. By corollary
\isom\ and the assumption just made, each $G_\epsilon$ is isomorphic
to $G_\a$ for some $\a\in I$.  By the pigeon hole principle there is  
a
single $G_\a$ which is isomorphic to uncountably many $G_\epsilon$.
This contradicts claim\NonIsom\endproof\main

\bigbreak
\neusection
\noindent
{\bf \S\number\secno\ Graphs without short cycles}

In this section we show that the class of all graphs omitting all
cycles of length at most $2k$ ($k\ge 2$) has no countable universal  
element.

\ppro \pentagon Definition: Let $S_k$ be the following graph: For  
five
vertices $\{x_i:i\in Z_5\}$ indexed cyclically connect $x_i$ to
$x_{i+1}$  by a simple path
$x_i,y_{i,1}\ldots,y_{i,k-1},x_{i+1}$.

PICTURE

\ppro \rigidThree Claim: If $f_1$ and $f_2$ are two embeddings of
$S_k$ into a graph $G$ omitting all  cycles of length at most $2k$  
($k\ge
2$) and $f_1(x_i)=f_2(x_i)$ for $i\in Z_5$ then $f_1=f_2$.

\proof: Suppose for simplicity that $f_1$ is the inclusion, and  
suppose
that $f_2\not=f_1$. Let $i$ be the least such that among  
$\{y_{i,j}:j<
k\}$ there is a vertex $v$ for which $v\not=f_2(v)$ and let $j(0)$ be
the least such that $y_{i,j}\not=f_2(y_{i,j(0)})$. Let $j(1)$ be the
the least $j>j(0)$ such that $x_{i,j(1)}=f_2(x_{i,j(1)}$. Now
$x_{i,j(0)-1},x_{i,j(0)},\ldots,x_{i,j(1)},f_2(x_{i,j(0)}),\ldots,
f_2(x_{i, j(1)-1})$ forms a cycle of length $\le 2k$ in $G$, contrary  
to
the assumption.\endproof\rigidThree

Let us define an infinite graph $U$ by induction. For every natural
$m$ let $S(m)$ be a copy of $S_k$ with vertices $x^m_{i}, y^m_{i,j}$  
($i\in
Z_5,1\le j < k$).

Let $U(0)=S(0)$. Suppose that $U(m)$ is defined
and  $S(m)\su U(m)$. To obtain $U(m+1)$ adjoin freely
 $S(m+1)$ of to $S(m)$ by identifying  $x^{m+1}_i$ with
 with  $y^m_{2i,[(k+1)/2]}$. Let $U=\bigcup U(m)$

Picture

\ppro \rigidFour Claim: (i) The graph $U$ contains no cycles of  
length
$\le 2k$, and $\deg_U(v)\le 3$ for all $v\in U$. (ii) If $f_1$ and  
$f_2$
are two embeddings of $U$ into a graph $G$ which contains no cycles  
of
length at most $2k$ and $f_1(x_i^0)=f_2(x_i^0)$ for $i\in Z_5$ then
$f_1=f_2$.

\proof: (i) is clear. Suppose that $f_1,f_2$ are as stated. Using
claim \rigidThree\ inductively one sees that  
$f_1=f_2$\endproof\rigidFour

Let us choose, by induction on $m$, vertices $v_m$ in $U$ such that  
the
distance in $U$ between $v_m$ and $v_{m+1}$ is at least $2k+1$.  For
every $\epsilon\in{}^\om 2$ let us construct a graph $U_\epsilon$ as
follows: let $u(m)$ be distinct vertices not in $U$. Connect $v(m)$  
to
$v(m+1)$ by an edge if $\epsilon(m)=1$ and connect $u(m)$ by edges to
$v(m),v(m+1)$ otherwise. Let $U_\epsilon=U\cup \{u(m):m\in M\}$.

It is not hard to verify that each $U_\epsilon$ omits all cycles of
length at most $2k$.

\ppro \cycleMain Theorem: For all $k\ge 2$ there is no universal
countable graph in the class of all graphs omitting all cycles of
length at most $2k$. In fact, the weak complexity of the class of all  
such
countable graphs is $\cont$.

\proof Suppose to the contrary that $\{G_\a:\a\in I\}$ is a set of
countable graphs, each omitting all cycles of length at most $2k$,
with the property that every countable graph omitting all cycles of
length at most $2k$ is isomorphic to a subgraph of $G_\a$ for at  
least
one $\a\in I$, and assume that $|I|<\cont$. Fix an embedding
$f_\epsilon$ of $U_\epsilon$ to some $G_{\a(\epsilon)}$. By the  
pigeon
hole principle there is a single $\a\in I$ which equals  
$\a(\epsilon)$
for all $\epsilon\in A$ for some uncountable set $A\su {}^\om 2$. For
each $\epsilon\in A$ let $\eta(\epsilon)=\lng f_\epsilon(x^0_i):i\in
Z_5\rng$. As there are only countably many finite sequences of
vertices in $G_\a$, there are different $\epsilon,\nu\in A$ with
$\eta(\epsilon)=\eta(\nu)$. But then it follows by claim \rigidFour\
that $f_\epsilon\rest U=f_\nu\rest U$. Let $m$ be such that
$\epsilon(m)\not=\nu(m)$. The vertices
$f_\epsilon(v(m)),f_\epsilon(u(m)),f_\epsilon(v(m+1))$ span a copy of
$C_3$ in $G_\a$, contrary to the assumption that $G_\a$ contains no
cycles of length at most $2k$. \endproof\cycleMain

\bigbreak
\bigbreak

\noindent{\bf References}

[CK] G.~Cherlin and P.~Komj\'ath, {\sl There is No Universal  
Countable
Pentagon-Free graph}, Journal of Graph Theory 18 (4) (1994) 337--341

[KP] Peter Komj\'ath and J\'anos Pach, {\sl Universal elements and  
the
complexity of certain classes of infinite graphs}, Discrete
Mathematics 95 (1991) 255--270

[KP1] P.~Komj\'ath and J.~Pach, {The complexity of a class of  
infinite
graphs}, Cominatorica 14(1) (1994) 121--125

[KMP] P.~Komj\'ath, A. Mekler and J. Pach, {\sl Some Universal
Graphs}, Israel Journal of Mathematics 64 (1988) 158--168

\end